\documentclass[runningheads]{llncs}
\usepackage[algo2e]{algorithm2e}
\usepackage{amsfonts}
\usepackage{amsmath}
\usepackage{amssymb}
\usepackage{caption}
\usepackage{dsfont}
\usepackage{graphicx}
\usepackage{hyperref}
\usepackage{nicefrac}
\usepackage{optidef}
\usepackage{subcaption}
\usepackage{xcolor}

\usepackage[export]{adjustbox}[2011/08/13]

\RestyleAlgo{ruled}

\newif\ifmargincomments %
\margincommentstrue

\ifmargincomments

\else

\fi

\newif\ifrelaxedv  %
\relaxedvfalse
\ifrelaxedv
\newcommand{\jvspace}[1]{}
 \else

\allowdisplaybreaks
\newcommand{\jvspace}[1]{\vspace{#1}}
\fi

\begin{document}
\title{Accelerating Continuous Variable \\ Coherent Ising Machines\\ via Momentum}
\author{Robin A. Brown\inst{1, 2, 3}\and
Davide Venturelli\inst{2, 3}\and \\
Marco Pavone\inst{1}\and
David E. Bernal Neira\inst{2,3,4}} 
\authorrunning{R. Brown et al.}
\institute{Stanford University, Autonomous Systems Laboratory \and
USRA Research Institute for Advanced Computer Science (RIACS) \and
NASA Quantum Artificial Intelligence Laboratory (QuAIL)  \and
Purdue University, Davidson School of Chemical Engineering}
\maketitle              %
\begin{abstract}
The Coherent Ising Machine (CIM) is a non-conventional architecture that takes inspiration from physical annealing processes to solve Ising problems heuristically.
Its dynamics are naturally continuous and described by a set of ordinary differential equations that have been proven to be useful for the optimization of continuous variables non-convex quadratic optimization problems.
The dynamics of such Continuous Variable CIMs (CV-CIM) encourage optimization via optical pulses whose amplitudes are determined by the negative gradient of the objective; however, standard gradient descent is known to be trapped by local minima and hampered by poor problem conditioning.
In this work, we propose to modify the CV-CIM dynamics using more sophisticated pulse injections based on tried-and-true optimization techniques such as momentum and Adam.
Through numerical experiments, we show that the momentum and Adam updates can significantly speed up the CV-CIM’s convergence and improve sample diversity over the original CV-CIM dynamics.
We also find that the Adam-CV-CIM’s performance is more stable as a function of feedback strength, especially on poorly conditioned instances, resulting in an algorithm that is more robust, reliable, and easily tunable.
More broadly, we identify the CIM dynamical framework as a fertile opportunity for exploring the intersection of classical optimization and modern analog computing.

\keywords{Ising Model \and Momentum \and Analog Computing.}
\end{abstract}
\section{Introduction}
Optimization is a powerful and intuitive framework for expressing and reasoning about the desirability of certain decisions subject to constraints.
With diverse applications ranging from supply chain management, vehicle routing, machine learning, and chip design, it is difficult to overstate the ubiquity of optimization.
In many contexts, casting a particular decision as an optimization problem--deriving mathematical expressions for the costs and constraints--permits the off-the-shelf application of commercial solvers as a satisfactory, general-purpose resolution to the problem.
However, even if a problem complies with a particular solver's specifications, practical constraints, such as memory or time, may bottleneck the size of problems that can be addressed or the quality of the solution obtainable.
While computational power has increased exponentially in the past few decades, according to Moore's law, if it stagnates these bottlenecks will persist, barring radical breakthroughs in computer engineering.

Quantum computing has emerged as an attractive and high-profile alternative to von Neumann computers based on its promise to accelerate specific computational tasks.
The past few years have seen significant commercial interest in quantum computing, leading to the rapid development of prototype quantum computers based on superconducting qubits \cite{kjaergaard2020superconducting}, photonics \cite{slussarenko2019photonic}, neutral atoms \cite{henriet2020quantum}, and trapped ions \cite{bruzewicz2019trapped}, encoding the gate-based model that underpins much of the theoretical work on quantum complexity theory.
However, due to the significant constant overhead induced by error correction, it is unlikely that quantum computing will result in practical speed-ups for optimization, notably when the best classical algorithms can be parallelized on specialized hardware, such as graphics processing units (GPUs) \cite{babbush2021focus}.
A competing paradigm to the gate-based model are quantum computers that can only solve particular forms of problems, or ``primitives", including quantum annealers (which solve Ising problems) \cite{hauke2020perspectives}, and Rydberg atom arrays (that encode the maximum independent set) \cite{ebadi2022quantum}.
While these are theoretically expressive primitives due to being NP-complete, the number of additional variables and connections required to construct a reduction from other NP-complete problems to one of these forms may significantly reduce the size of problems addressable on near-term hardware, despite only being polynomial overheads \cite{lucas2014ising}, \cite{kleinbergtardos2005}.
It is against this backdrop that there has been a recent surge of interest in non-quantum, yet non-conventional computing architectures for optimization based on coupled lasers \cite{pal2020rapid}, memristors \cite{cai2020power}, polaritons \cite{berloff2017realizing} \cite{kalinin2020polaritonic}, and optical parametric oscillators \cite{marandi2014network} \cite{mcmahon2016fully}. 

From an algorithmic perspective, the development of these architectures has concurrently ushered in research on hybrid algorithms that explore how non-conventional and classical computers should complement each other \cite{callison2022hybrid}.
These algorithms are typically developed with specific hardware abstractions in mind, so their practical impact will depend heavily on which technologies prevail in a highly dynamic field.
In this paper, we pose a complementary question of how classical optimization algorithms can inform the development of non-conventional architectures.
Despite their theoretical limitations, classical optimization techniques have achieved impressive empirical success, even for difficult, non-convex optimization problems.
Rather than discarding the innovations developed in the context of classical optimization, we contend that these ideas should be incorporated into future computing architectures.
Crucially, much of the discussion of these technologies has been couched in the language of physics and optics, making it largely inscrutable for the optimization community and obscuring the opportunities to incorporate well-known optimization techniques.

This paper focuses on one such device, the Coherent Ising Machine (CIM), which takes inspiration from the principles of thermalization to search for the ground state of Ising problems heuristically.
It was first prototyped in 2014 as a network of four optical parametric oscillators (OPOs) coupled by optical delay lines \cite{marandi2014network}.
Since then, the technology has rapidly developed, with prototype devices reaching sizes of 100,000 spins \cite{honjo2021100}.
One significant deviation in today's architecture versus the original architecture is the inclusion of a measurement-feedback circuit implemented with a field programmable gate array (FPGA) \cite{mcmahon2016fully}.
Not only has this development enabled scalable all-to-all connections between the OPOs, it has also enabled the implementation of more sophisticated feedback methods, which have primarily been focused on mitigating the ``amplitude heterogeneity" problem observed in early instantiations of the device \cite{leleu2017combinatorial}.
This paper focuses on the CIM because the measurement-feedback circuit allows for modular modifications to the feedback term without otherwise affecting the operational principles of the device, making it a fertile opportunity for exploring the intersection of classical optimization and analog computing.

\emph{Contribution:}
Our contribution hinges on the insight that all existing variations of the CIM rely on gradient descent to encourage the optimization of a desired objective.
In this paper, we contend that the feedback term should be modularly swapped with more sophisticated optimization techniques, particularly the accelerated and adaptive optimizers that are tried and true for large-scale non-convex optimization, such as neural network training.
We conduct an extensive numerical evaluation of the proposed modifications to evaluate the anticipated performance gains realizable in a physical device. 

\emph{Organization:} 
In Section \ref{sec:background}, we will first define our notation and terminology.
We will then present a high-level overview of the CIM's operational principles, along with relevant developments, and conclude the section with a discussion of the momentum and Adam optimizers.
In Section \ref{sec:dynamics}, we will present the Gaussian-state model of the CIM's dynamics, highlighting opportunities to incorporate other optimization techniques. 
We will also discuss the specific application of the CIM to box-constrained quadratic programs, which we later use for benchmarking.
In Section \ref{sec:experiments}, we present numerical experiments evaluating convergence speed, diversity of solutions, and parameter sensitivity of the modified CIM against its original dynamics.

\section{Background and Related Work}\label{sec:background}
In this section, we will first present our notation and terminology.
We will then provide a high-level overview of the CIM's operational principles, discussing its history and recent developments designed to improve optimization performance.
Finally, we will briefly introduce momentum and the Adam optimizer, overviewing some of their desirable properties, particularly for non-convex optimization.

\subsection{Notation and Terminology}

In this paper, we work with vectors and matrices defined over the real numbers and reserve lowercase letters for vectors and uppercase letters for matrices.
We will also follow the convention that a vector $x \in \mathbb{R}^n$ is to be treated as a column vector, i.e., equivalent to a matrix of dimension $n \times 1$.
For a matrix $M$, we use $M_{i, j}$ to denote the entry in the $i$th row and $j$th column, $M_{i, *}$ denotes the entire $i$th row, and $M_{*, j}$ denotes the entire $j$th column.
An Ising problem is an optimization problem of the form: $\min_{x \in \{-1, 1\}^n} f(x) := \sum_{i, j} Q_{i, j} x_i x_j + \sum_{i} c_{i} x_i$, where $Q \in \mathbb{R}^{n \times n}, \, c \in \mathbb{R}^n $ are real coefficients and $x_i \in \{-1, 1\}$ are discrete variables to be optimized over.
We will use the convention that the discrete variables are $\pm 1$, also known as \emph{spin} variables.

\subsection{Coherent Ising Machines}
The Coherent Ising Machine (CIM) is an optical network of optical parametric oscillators (OPOs) that was originally designed to heuristically search for the solution (``ground state") of Ising problems.
At a high level, the CIM relaxes the discrete variables to continuous variables and augments the Ising objective, $f(x)$, with a ``double-well potential" of the form $\phi (x, a) := \sum_i \frac{1}{4}x_i^4 - \frac{a}{2}x_i^2$. 
The OPOs, whose amplitudes represent the variables of the Ising problem, are encouraged to optimize the Ising objective by evolving approximately in the negative gradient direction of the augmented objective.
At the same time, the parameter $a$ (the ``pump term”) is gradually increased (``annealed") during the device's run.
This has the effect of deepening the wells of the potential and encouraging the OPOs to take on the same magnitude.
While this closely resembles the standard penalty method (i.e., approximating constraints with a penalty term in the objective), it also has the effect of gradually reshaping the energy landscape from a trivial form to one that encodes the desired objective. 
Moreover, the wells of the potential are at $\pm \sqrt{a}$, so the mapping from OPO amplitude to binary variables also varies while the pump term is annealed. 

The double-well potential is physically implemented using a pumped laser and a periodically-poled lithium niobate (PPLN) crystal \cite{marandi2014network}. 
Early prototypes encoded the objective with optical delay lines and only accommodated objectives with linear gradients (i.e., quadratic functions), allowing the CIM to naturally encode Ising problems, as its name suggests.
Despite the simplicity of their functional form, Ising problems are a rich class of problems encompassing Karp's NP-complete problems \cite{lucas2014ising}, and NP-hard problems such as job-shop scheduling \cite{venturelli2016job}, vehicle routing \cite{harwood2021formulating}, and community detection \cite{negre2020detecting}, to name a few.

While significant progress has recently been made towards deriving a theory of the (fully classical) CIM's performance on Sherington-Kirkpatrick spin glass problems \cite{yamamura2023geometric}, a complete theory characterizing its performance on generic optimization problems still remains elusive.
Currently, the CIM is primarily benchmarked by simulating a quantitative model of its behavior on different applications, such as MIMO detection \cite{kim2021physics} and compressed sensing \cite{gunathilaka2023effective}.
Even though this is a widely accepted approach, no single model of the CIMs dynamics exists.
Instead, different models have been constructed with varying degrees of fidelity when modeling the quantum-mechanical effects.

In parallel with the theoretical developments, the hardware of the CIM has matured considerably from the early prototypes \cite{honjo2021100}.
Instead of implementing couplings between spins using optical delay lines, couplings are now implemented using a measurement-feedback circuit designed to allow for scalable all-to-all connections between OPOs \cite{mcmahon2016fully}.
This circuit consists of a measurement process (optical homodyne detection) to determine the amplitudes of the OPOs, and a feedback term based on these measurements computed on a field-programmable gate array (FPGA).
In this work, we propose to leverage the flexibility afforded to us by the measurement-feedback circuit to modify the feedback computation.

While modifying the feedback term has been explored previously, efforts have primarily been targeted toward heuristic resolutions to the ``amplitude heterogeneity" problem.
This occurs when the OPO amplitudes are not exactly $\pm \sqrt{a}$, breaking the mapping between the CIM and Ising objectives.
From an optimization perspective, this can be seen as the infeasibility that may occur when replacing a constraint with a penalty function.
One such modification is the amplitude heterogeneity correction proposed in \cite{leleu2019destabilization}, which introduced an auxiliary error correction variable that exponentially increases/decreases if the spin amplitudes are greater/less than $\sqrt{a}$.
These error variables modulate the relative weights of the Ising objective and double-well potential and increase the latter when the spins have different amplitudes.
Later work extended this idea to additionally modulate the pump term, $a$, to encourage the CIM to continue exploring the solution space within a single run rather than converging to a single solution \cite{kako2020coherent}.
A fully optical implementation of the error correction mechanism was later developed in \cite{reifenstein2021coherent}, thus overcoming some of the downsides of having an electronic component in the loop.

The focus on resolving the amplitude-heterogeneity problem has primarily been motivated by problems whose complexity stems from optimizing over discrete domains.
This starkly contrasts the many applied works that artificially discretize continuous variables to force them into the Ising form \cite{borle2019analyzing}\cite{chang2019quantum}\cite{ottaviani2018low}\cite{willsch2020support}. 
This encoding dramatically increases the number of OPOs mapped to a single variable and may result in a significant skew in the coupling coefficients.
Recognizing the expressive power of \emph{continuous} optimization problems, \cite{khosravi2022non} recently proposed embracing the naturally continuous nature of the OPOs to solve box-constrained quadratic programs (BoxQPs), whose complexity stems from non-convexity rather than discreteness--this modified device is called the coherent continuous variable machine (CCVM).
The CCVM does not fundamentally change the operational principles of the CIM but instead interprets the OPO amplitudes as their native continuous values rather than projecting them to binary values.
While this primarily impacts how the resulting solution is interpreted, it can materially affect the dynamics of the device, as the feedback term is also computed with the continuous values rather than their binary projection. 
This modification was further supported by \cite{brown2022copositive}, who showed that BoxQPs are an algorithmic primitive that enables copositive optimization, a broad class of optimization problems that includes Ising problems.

\subsection{Accelerated and adaptive optimizers}
Given an optimization problem, $\min_x f(x)$, gradient descent is one of the most primitive and fundamental algorithms one could use to solve it.
It formalizes the intuition that the optimization variable should be iteratively updated in the direction of the greatest decrease in the objective.
The variables are updated according to $x^{(t + 1)} = x^{(t)} - \gamma_t \nabla f(x^{(t)})$, where $\gamma_t > 0$ are iteration-specific step-sizes.
Stochastic gradient methods replace $\nabla f(x^{(t)})$ with a cheap, unbiased estimate, and are among the most popular methods for large-scale optimization due to their broad applicability, simplicity, and efficiency.

One scenario known to cause slow convergence for gradient descent is when the objective function is poorly-conditioned. 
Intuitively, the conditioning of an optimization problem refers to the disparity in sensitivity of the objective function along different directions.
A well-conditioned problem is one where the objective function has roughly the same sensitivity in all directions.
In contrast, in a poorly conditioned problem, the objective is far more sensitive along some directions than others.
In the latter case, gradient descent exhibits oscillatory behavior along the directions where the objective is most sensitive while making little progress along directions of low sensitivity. 

Momentum methods, named for the analogy to mechanical systems, have been shown to improve convergence for poorly conditioned problems by updating the variable with an exponentially decreasing weighted sum of all past gradients.
Mathematically, the updates are of the form
\begin{align}
    g^{(t + 1)} &= \theta g^{(t)} + (1 - \theta) \nabla f(x^{(t)}), \label{eq:mom_update}\\
    x^{(t + 1)} &= x^{(t)} - \gamma g^{(t + 1)}, 
\end{align}
where $\theta \in [0, 1)$ is a damping parameter.
In poorly conditioned problems, this has the effect of smoothing out oscillations in the high sensitivity direction while adding up contributions along the low sensitivity direction, thereby improving convergence speed \cite{qian1999momentum}.
It is also effective for helping iterates escape saddle-points and local minima, which can trap ordinary gradient descent methods.

In principle, if poor conditioning is a result of inconsistent scaling between variables, it could be combated by using different step sizes per variable.
However, this introduces additional challenges, as the step sizes must be individually chosen for each variable.
Instead, the Adam optimizer uses estimates of the first and second moments of the gradient to adapt the step size for each variable automatically.
The complete algorithm is presented in Algorithm \ref{alg:adam}.
Adam is well-known for its empirical success on large-scale, non-convex optimization problems, such as neural network training, for which it is often chosen as the default optimizer.

\begin{algorithm2e}
\caption{Adam optimizer updates}\label{alg:adam}
\KwIn{$f(\cdot)$ function to be optimized; Parameters $\beta_1, \beta_2 \in [0, 1)$;\\ Initial variable vector $x^{(0)}$}
\KwOut{$x^{(t)}$ that approximately minimizes $f(x)$}
$t \gets 0, \quad m^{(0)} \gets 0, \quad v^{(0)} \gets 0$ \;
\While{$x^{(t)}$ \text{has not converged}}{
    $g^{(t)} \gets \nabla_x f(x^{(t)})$\;
    $m^{(t)} \gets \beta_1 m^{(t - 1)} + (1 - \beta_1) g^{(t)}$\;
    $v^{(t)} \gets \beta_2 v^{(t - 1)} + (1 - \beta_2) (g^{(t)})^2$\;
    $\hat{m}^{(t)} \gets \frac{m^{(t)}}{1 - \beta_1^t}$\;
    $\hat{v}^{(t)} \gets \frac{v^{(t)}}{1 - \beta_2^t}$\;
    $x^{(t)} \gets x^{(t - 1)} -  \gamma \frac{\hat{m}^{(t)} }{\sqrt{\hat{v}^{(t)}} + \epsilon}$\;
    $t \gets t + 1$
}
\Return{$x^{(t)}$}\;
\end{algorithm2e}

\section{Dynamics of the CIM}\label{sec:dynamics}
In this section, we will present the Gaussian-state model of the CIM's dynamics, which will be the basis of our numerical experiments.
While this is not a fully quantum treatment of the device, as it only considers up to second-order quantum correlations, it has been found to be consistent with experimental CIMs and models important phenomena such as squeezing/anti-squeezing, measurement uncertainty, and backaction, while maintaining the computational tractability necessary to simulate large-scale systems \cite{ng2022efficient}.
Under this model, it is assumed that each OPO pulse amplitude is a Gaussian distribution with mean $\mu_i$ and variance $\sigma_i$ that evolves according to the following dynamics \cite{kako2020coherent}:
\begin{subequations}\label{eq:cim_system}
\small{
    \begin{equation}\label{eq:mu_update}
    \begin{aligned}
        \mu(t + \Delta t) = \mu(t) & + \bigg[(-(1 + j_t) + p_t - g^2 \mu(t)^2)\mu(t) 
                            - \lambda \partial f(\tilde{\mu}(t))\bigg]\Delta t \\
                            & +\left[ \sqrt{j_t}(\sigma(t) - 1/2)\right] Z(t) \sqrt{\Delta t},\\
    \end{aligned}
    \end{equation}
    \begin{equation}\label{eq:sigma_update}
        \begin{aligned}
        \sigma(t + \Delta t) = \sigma(t) &+ 
        \bigg[2(-(1 + j_t) + p_t - 3g^2\mu(t)^2)\sigma(t) 
        -2j_t(\sigma(t) - 1/2)^2\\
        &+ ((1 + j_t) + 2g^2\mu(t)^2)\bigg] \Delta t,
        \end{aligned}
    \end{equation}
    \begin{equation}\label{eq:mu_meas}
    Z_i(t) \sim \mathcal{N}(0, 1), \qquad
    \tilde{\mu}(t) = \mu(t) + \frac{Z(t)}{\sqrt{4j_t\Delta t}}.
    \end{equation}}
\end{subequations}
The remaining quantities in the above equations represent the normalized continuous measurement strength, $j_t$,  a normalized non-linearity coefficient, $g$, the pump strength, $p_t$, and feedback magnitude, $\lambda$ (which is analogous to the step-size in optimization algorithms).
While these quantities are governed by the physical parameters of the device, in this work, we treat them abstractly as hyperparameters that can be optimized to improve the performance of the CIM.
Following the treatment in \cite{khosravi2022non}, the hyperparameters $j_t$ and $p_t$ are indexed by $t$ to emphasize that they are modulated over the run of the device.

In Equation \eqref{eq:mu_update}, the term $(-(1 + j_t) + p_t - g^2 \mu(t)^2)\mu(t)$ is the gradient of a double-well potential with minima at $\pm \nicefrac{\sqrt{ p_t  - (1 + j_t)}}{g}$.
At each iteration (``round trip"), the mean-field amplitude is measured using optical homodyne detection. 
The measured amplitude, presented in Equation \eqref{eq:mu_meas}, may differ from its true value by a finite uncertainty depending on the measurement strength.
The measurement also induces a shift in the mean-field amplitude, represented by the term $\sqrt{j_t}(\sigma(t) - 1/2)Z(t) \sqrt{\Delta t}$, through backaction.
The feedback term, $-\lambda \partial f(\tilde{\mu}(t))$, is computed based on the measured amplitudes.
Crucially, for measurement-feedback variants of the CIM, this term is computed on a field programmable gate array (FPGA) rather than being implemented in optics and thus can be modified modularly without otherwise altering the hardware.

Abstractly, this is represented by replacing the term $\partial f(\tilde{\mu}(t))$ in {Equation~\eqref{eq:mu_update}} with an arbitrary update $\Phi(\{\tilde{\mu}(t) \}_t)$.
Because updates that are too complex or compute-intensive may undermine the resource bottlenecks motivating the CIM in the first place, we aim to design $\Phi$ to enhance the optimization performance of the system defined by Equation \eqref{eq:cim_system} while maintaining comparable complexity to the feedback computations in the amplitude-heterogeneity-correction (AHC) and chaotic-amplitude-control (CAC) variants of the CIM.
In both variants, the complexity of the update is dominated by the complexity of computing the gradients, $O(n^3)$.
At the same time, additional computations induce an additional overhead of $O(n)$, where $n$ is the dimension of the decision variable $x$.
Specifically, the momentum and Adam updates from Equation \eqref{eq:mom_update} and Algorithm \ref{alg:adam}, respectively, satisfy these desiderata.
In the remainder of the paper, we will refer to the CIM with momentum or Adam updates as momentum-CV-CIM and Adam-CV-CIM, respectively.
The original CIM variant will be referred to as GD-CV-CIM to emphasize that the updates are based on gradient descent.

\subsection{Application to box-constrained quadratic programs}
In this work, we will benchmark the proposed modifications of the continuous CIM on box-constrained quadratic programs (BoxQP).
These are optimization problems of the form:

\begin{mini}
{x \in \mathbb{R}^n}{x^\top Q x + c^\top x}
{\label{eq:boxqp}\tag{BoxQP}}{}
\addConstraint{\textbf{0} \leq x \leq \mathds{1},}
\end{mini}

where there are no restrictions on $Q \in \mathbb{R}^{n \times n}$ or $c \in \mathbb{R}^n$.
\eqref{eq:boxqp} has been shown to be NP-hard if $Q$ is negative semi-definite or indefinite \cite{pardalos1991quadratic}.
Furthermore, \eqref{eq:boxqp} is an expressive primitive that enables solving a broader class of mixed-binary quadratic programs, including Ising problems \cite{brown2022copositive}.
In this work, we seek to evaluate whether momentum or Adam updates improve the convergence speed of the CIM dynamics, reducing the time needed to obtain each sample from the device.
Evaluating the CIM based on the continuous values of the OPO amplitudes will allow us to detect mechanistic improvements to convergence without the confounding and discontinuous process of projecting OPO amplitudes to binary variables based on their sign.

While the double-well potential induces a bias for the amplitudes $\pm \sqrt{a}$ (in Equation \eqref{eq:mu_update},  $a = \nicefrac{p_T - (1 + j_T)}{g^2}$) we will map the entire domain of the OPOs, $[-\sqrt{a}, \sqrt{a}]^n$, to the domain of \eqref{eq:boxqp}.
This is given by the transformation, 
\begin{equation}
\small{
    x_i = \max \left(
    \min \left(\frac{1}{2}\left(\frac{\mu_i}{\sqrt{a}} + 1\right), 
    1 \right),
    0\right).}
\end{equation}
Intuitively, this re-scales and shifts the range $[-\sqrt{a}, \sqrt{a}]$ to $[0, 1]$ and projects any infeasible variables to the feasible domain.
We note that the double-well potential only softly confines the OPOs to this domain (akin to penalty functions), so the OPOs may exceed the feasible region.
It also introduces an extraneous quadratic term due to its inception as an Ising solver rather than a continuous optimization solver.
While we expect this to degrade the performance of the resulting solver, our objective is to assess the proposed modifications on already realized instantiations of the CIM and thus should be reflected by the model.

\section{Experiments}\label{sec:experiments}
In Section \ref{sec:dynamics}, we presented a model of the CIM's dynamics, highlighting an opportunity to incorporate more sophisticated optimization techniques.
In this section, we aim to assess whether such substitutions have any performance benefit for the CIM.
The performance of quantum/quantum-inspired devices is often evaluated on the time-to-target metric \cite{ronnow2014defining}.
The time-to-target metric with $s$ confidence represents the number of repetitions to sample a state that achieves a particular ``target" with probability $s$, multiplied by the time for a single run of the algorithm, $\texttt{T}_{\text{single}}$, i.e., 
\begin{align}\label{eq:ttt}
\small{
    \texttt{TTT}_s = \texttt{T}_{\text{single}}\frac{\log(1 - s)}{\log(1 - \hat{p}_{\text{succ}})}.
}
\end{align}
In the context of continuous optimization, the targets typically correspond to achieving a particular relative optimality gap, i.e.,  $\text{gap}_f(x) := \frac{f(x) - f(x^*)}{\vert f(x^*) \vert} \leq \epsilon$.

The time-to-target metric incorporates two figures of merit: the time to convergence of a single run of the algorithm and the distribution to which the samples from the algorithm converge.
If a particular method improves this metric, it can be due to any combination of these two mechanisms.
First, the samples can converge to their final state more quickly, resulting in a decrease in the multiplier $\texttt{T}_{\text{single}}$.
Secondly, the method may change the distribution that the samples converge to, effectively changing the probability of success $\hat{p}_{\text{succ}}$.
Our numerical experiments are intended to evaluate the performance of the modified CIM against the original variant along these two metrics.
In addition, we will evaluate the robustness of the different variants as a function of the feedback strength.
This is one component of the effort required for a practitioner to successfully apply the CIM and its variants on an unseen problem instance; this is an important, yet often neglected, factor in the practical applicability of a solver.

In this section, we report results on the standard benchmark instances from \cite{Burer10},\cite{BurVan09},\cite{VanNem05b}, available at \href{https://github.com/sburer/BoxQP_instances}{https://github.com/sburer/BoxQP\_instances}.
All code required to replicate our experiments is available at \href{https://github.com/SECQUOIA/continuous-CIM}{https://github.com/SECQUOIA/continuous-CIM}.

\subsection{Time to optimality gaps}\label{subsec:timefrac}
\begin{figure}
    \centering
    \includegraphics[width=0.5\textwidth]{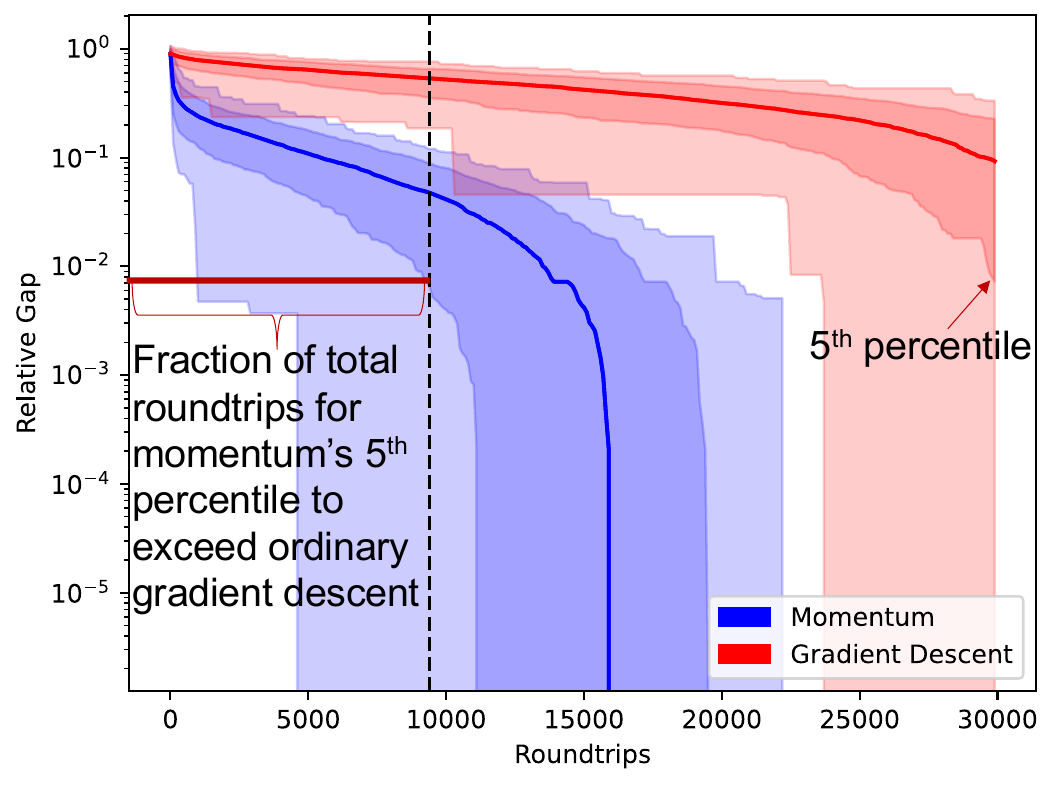}
    \caption{
    This figure shows how the roundtrip ratio is computed for a single problem instance.
    We determine the best relative optimality gap achieved by the 5th percentile of both CIM variants (8e-3 at the end of GD-CV-CIM's evolution) and compare the roundtrip numbers where the 5th percentile first reaches this gap ($\sim$ 9500 for momentum, 30000 for gradient descent).
    The roundtrip ratio is approximately 0.32, and we indicate which variant achieves the target gap faster.}
    \label{fig:time_frac_ex}
\end{figure}

The first experiment evaluates the relative convergence speed of the modified CIMs versus the original dynamics.
For each problem instance, we evaluate the relative optimality gaps for the $X$th percentile, with $X = [5, 10, 25, 50]$, throughout the CIM's evolution and determine the best gap achieved by the $X$th percentile of both variants (hereafter referred to as the target gap).
We then evaluate the roundtrip where each variant achieves the target gap for the first time.
We record the ratio of the faster variant's roundtrip value to that of the slower variant, along with which one was faster.
Figure \ref{fig:time_frac_ex} plots how this metric is calculated for a single problem instance.
Intuitively, the roundtrip fraction represents the relative improvement in $\texttt{T}_{\text{single}}$.

We note that the roundtrip fraction metric does not capture whether a sample achieves its minimum optimality gap only transiently or if the dynamics stabilize around that value.
Consequently, this metric makes the most sense in contexts where the optimality gaps are close to decreasing monotonically.
In this paper, we leave the OPO bounds to their ``natural" values of $[-\sqrt{a}, \sqrt{a}]$, with $a = \frac{p_T - (1 + j_T)}{g^2}$.\footnote{The remaining parameters are given by:\\ $T = 15000, \, dt=0.0025, \, g = 0.01, \, \lambda = 550, \, p_t = 2.5 \frac{t}{T}, \, j_t = 25 \exp{\left(-3 \frac{t}{T}\right)}$.}
Setting the bounds to be significantly smaller than those encouraged by the double-well potential means that the OPO amplitudes are primarily evolving outside the feasible domain.
This induces fast transients in optimality gaps when a particular variable can be updated from one side of the feasible domain to another within a single iteration.  
Because the roundtrip fraction metric is designed to measure differences in mechanistic convergence to a solution, the bounds are left to the natural domain of the OPOs, instead of the aggressive clipping proposed in \cite{khosravi2022non}.

\begin{figure}
    \centering
    \includegraphics[width=0.8\linewidth]{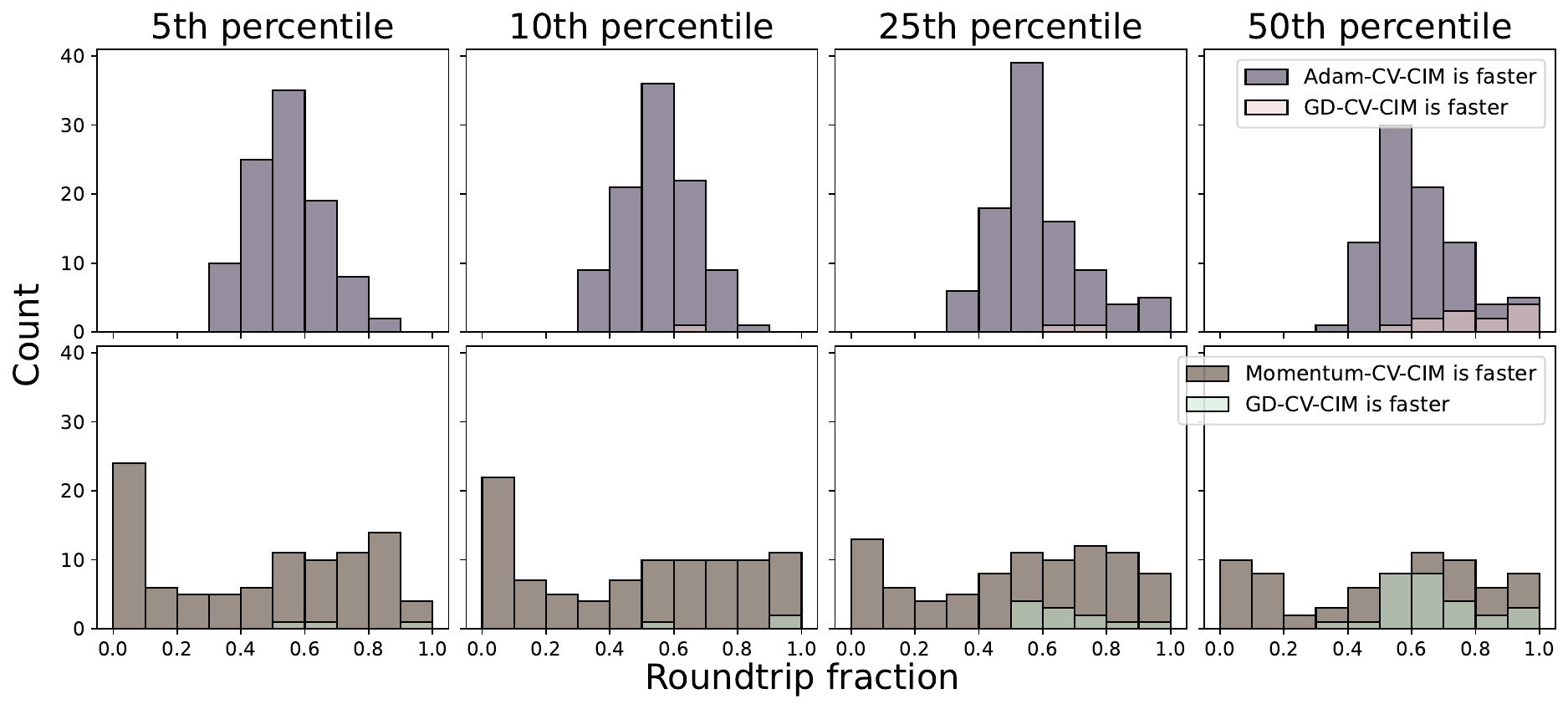}
    \caption{This figure plots the histogram of roundtrip fractions, computed according to the process illustrated in Figure \ref{fig:time_frac_ex}.
    The top row compares the Adam-CV-CIM to the GD-CV-CIM, while the bottom row compares the momentum-CV-CIM to the GD-CV-CIM.
    The data points are categorized depending on which CIM feedback variant converged faster.}
    \label{fig:timefrac_single}
\end{figure}

Figure \ref{fig:timefrac_single} plots the histogram of roundtrip fractions, with the top row comparing the Adam-CV-CIM to the GD-CV-CIM and the bottom row comparing the momentum-CV-CIM to the GD-CV-CIM.
We find that both the momentum-CV-CIM and Adam-CV-CIM consistently reach the target gap before the GD-CV-CIM does. 
However, for Adam-CV-CIM, the roundtrip fractions are normally distributed around 0.5 (meaning the Adam-CV-CIM is $\sim2$ times faster), while for the momentum-CV-CIM, the roundtrip fractions are more uniformly distributed across the entire range of values with a slight peak in the 0 -- 0.1 bucket ($\geq$10 times faster than the GD-CV-CIM).
The difference between the momentum-CV-CIM and ordinary CIM is most pronounced for the top percentiles of samples, while it is less pronounced for the 50th percentile (i.e., median of samples).
In contrast, the relative performance of the Adam-CV-CIM is more consistent across percentiles.

\subsection{Sample distribution}\label{subsec:distribution}
\begin{figure}
    \centering
    \includegraphics[width=0.7\linewidth]{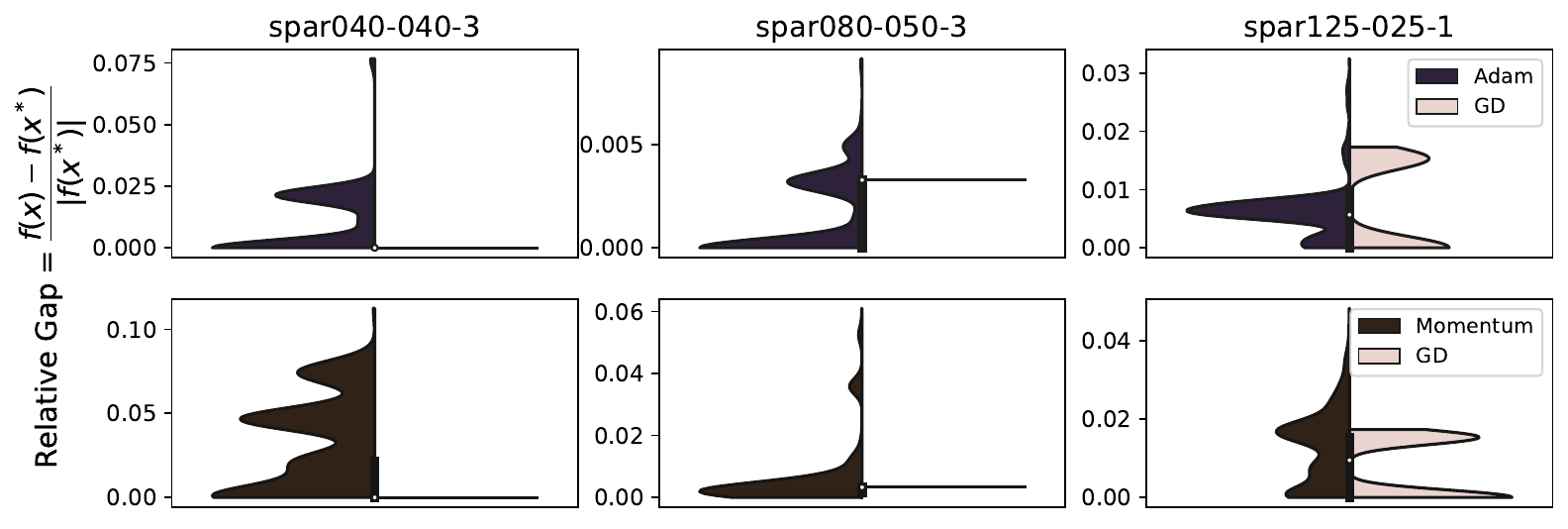}
    \caption{This figure plots the relative optimality gaps for three representative instances.
    In all instances, the Adam-CV-CIM and momentum-CV-CIM achieve a greater diversity of samples with multi-modal distributions.
    For \texttt{spar040-040-3}, all samples from GD-CV-CIM concentrate on the global optimum, while for \texttt{spar080-050-3}, all samples concentrate on a suboptimal solution.
    For \texttt{spar125-025-1}, all variants result in multi-modal distributions with some shared modes between solvers.}
    \label{fig:gap_violin}
\end{figure}

Secondly, we sought to understand whether there was a difference in the distributions found by the modified CIMs.
Figure \ref{fig:gap_violin} depicts violin plots of the relative optimality gaps for three representative instances.
In all instances, we find that the Adam-CV-CIM and momentum-CV-CIM achieve a greater diversity of samples with multi-modal distributions.
In contrast, for the \texttt{spar040-040-3} and \texttt{spar080-050-3} instances, all samples from the GD-CV-CIM achieve the same optimality gap; they concentrate on the global optimum for \texttt{spar040-040-3}, but only achieve a relative optimality gap of 0.003275 for \texttt{spar080-050-3}.
For the instance \texttt{spar125-025-1}, all variants result in multi-modal distributions, which share peaks in optimality gaps. 
The momentum-CV-CIM has three peaks in its distribution with relative optimality gaps of 0, 0.007, and 0.018.
The GD-CV-CIM shares two of these peaks at 0 and 0.018, while the momentum-CV-CIM shares the peaks at 0 and 0.007.

\begin{figure}[ht!]
    \centering
    \includegraphics[width=0.8\linewidth]{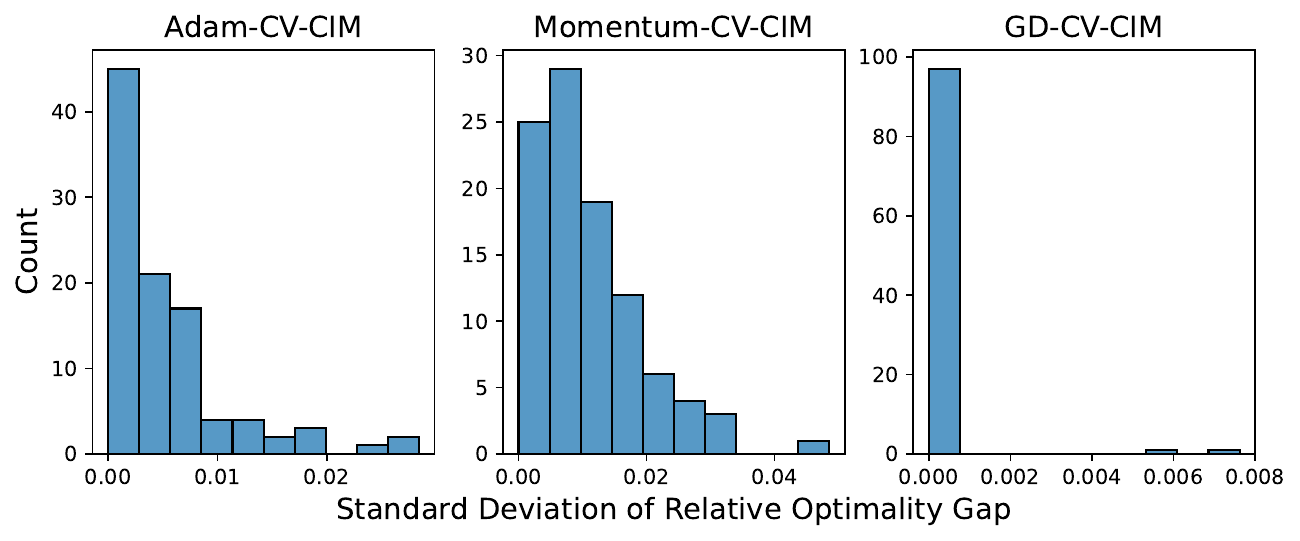}
    \caption{This figure plots the histogram of optimality gap standard deviations for all CIM variants.
    The momentum-CV-CIM has the highest standard deviations, with values roughly double that of the Adam-CV-CIM.
    The GD-CV-CIM has very low variation in optimality gaps, with most instances having no variation.}
    \label{fig:gap_variances}
\end{figure}

We aggregate these metrics over all problem instances by computing the standard deviation in the relative optimality gap achieved by each sample.
Figure \ref{fig:gap_variances} plots the histogram of optimality gap standard deviations for all CIM variants.
We find that the momentum-CV-CIM has the highest standard deviations, with values roughly double that of the Adam-CV-CIM.
In contrast, consistent with the examples depicted in Figure \ref{fig:gap_violin}, the GD-CV-CIM has very low variation in optimality gaps, with most instances having no variation.
We note that the relationship between the degree of variation and solver performance is not monotonic.
In the context of Equation \eqref{eq:ttt}, if all samples converge to the same value, $\hat{p}_{\text{succ}}$ will be either 0 or 1, depending on the desired target.
This means that the solver will either achieve the target with a single sample or not at all.
While this may improve the time-to-target for problems that are particularly amenable to the solver, it will result in a complete failure for problems that are not.
This is contrary to the perspective of the CIM as a sampler, from which one may expect a diverse sample of high-quality, albeit sub-optimal, solutions.
Instead, the solver should strike a delicate balance between achieving a good diversity of samples while still concentrating on those with good objective value.

\begin{figure}[ht!]
    \centering
    \includegraphics[width=0.8\linewidth]{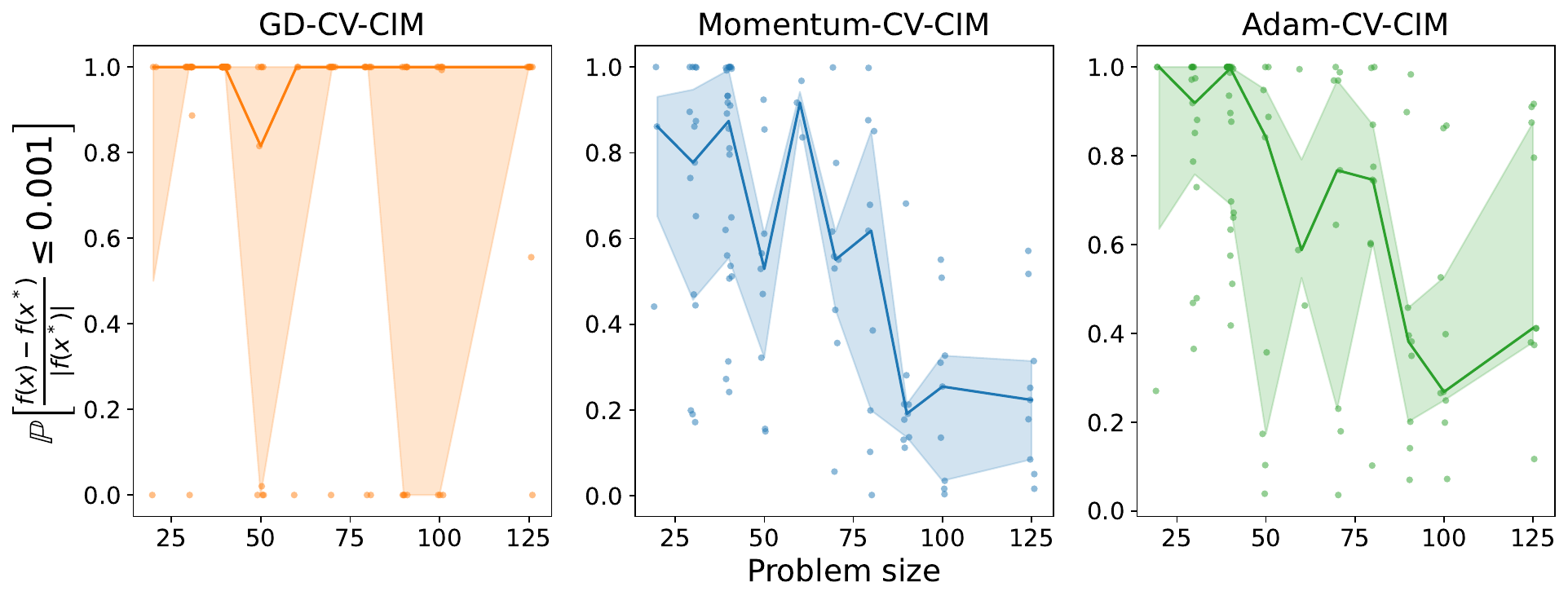}
    \caption{This plot shows the likelihood of samples having a relative optimality gap $\leq 0.1 \%$. Adam-CV-CIM and momentum-CV-CIM consistently achieved this, while GD-CV-CIM often met the target, except in 17 instances out of 99.}
    \label{fig:boxqp_psucc}
\end{figure}

Figure \ref{fig:boxqp_psucc} plots the probability that each variant finds samples with a relative optimality gap of at most $0.1\%$.
On the standard benchmarks, using default parameters, both the Adam-CV-CIM and momentum-CV-CIM had non-zero probabilities of finding a solution within a $0.1\%$ gap on all instances.
In contrast, the GD-CV-CIM frequently had all samples converging within a $0.1\%$ gap, however in 17 (out of 99 total) instances, the GD-CV-CIM failed to find any solutions within a $0.1\%$ gap.

\subsection{Sensitivity to feedback strength}
In Sections \ref{subsec:timefrac} and \ref{subsec:distribution}, we found that all variants of the CIM had reasonable performance on the standard benchmarks, even using default parameters.
However, one commonly neglected factor is the degree to which solver performance depends on its hyperparameters.
Indeed, performance metrics are often reported based on hyperparameters that have been meticulously hand-tuned, typically in an opaque process that obscures the effort required for a practitioner to apply the solver to a new problem.
In this section, we introduce a new benchmark set designed to test the CIM variants on problems with poor conditioning.
At a high level, these instances are generated to have a non-convex landscape and are additionally parameterized by a constant $\kappa \in \mathbb{R}_{\scriptsize{++}}$, that dictates the skew when re-scaling the objective along different axes.\footnote{Specifically, instances are constructed as $M = D(\kappa) U \Sigma(n) U^\top D(\kappa)$, where $U \in \mathbb{R}^{n \times n}$ is a random orthogonal matrix, $\Sigma(n) \in \mathbb{R}^{n\times n}$ is a diagonal matrix with $\Sigma(n)_{ii} = 1$  if $i \leq \frac{n}{2}$ and $-1$ otherwise (to induce a non-convex landscape), and ${D(\kappa) = \text{diag}(1, 1 + \frac{\kappa - 1}{n} , \ldots , \kappa)}$.}

\begin{figure}[ht!]
    \includegraphics[width=0.9\textwidth, center]{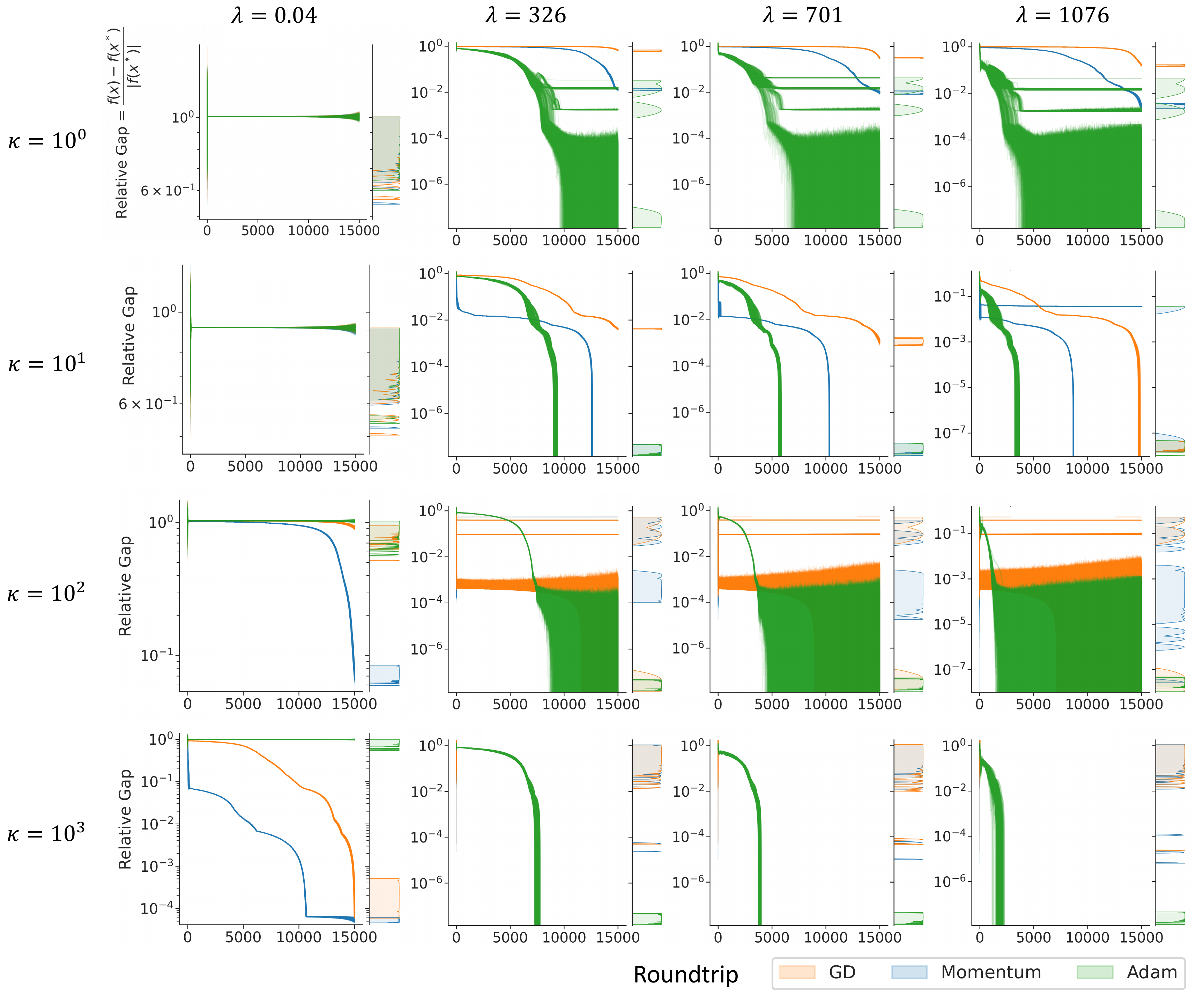}
    \caption{Each plot in the grid illustrates the convergence in relative optimality gaps over roundtrips for varying combinations of $\kappa$ and feedback strengths, $\lambda$. 
    We plot the kernel-density estimate of the best optimality gaps on the right axis of each subplot to illustrate their distributions.
    Notably, the GD and momentum CV-CIMs show limited consistency in successful feedback strengths with respect to problem size, while the Adam-CV-CIM exhibits more consistent dynamics across instances of varying $\kappa$.}
    \label{fig:stepsize}
\end{figure}

In this section, we evaluate the sensitivity of each of the variants to the feedback strength, $\lambda$ in Equation \eqref{eq:mu_update}, analogous to step-size in classical optimization.
Figure \ref{fig:stepsize} plots the convergence in relative optimality gaps for representative instances of dimension 20 with $\kappa \in [1, 10, 1000, 1000]$ across feedback strengths of $\lambda \in [0.04, 326, 701, 1076]$.
We plot the kernel-density estimate of the best gap per sample on the right axes. 
We note that absence of samples from an individual plot indicates that the sample diverged, and the best gap is computed until the point of divergence.
We find that there is limited consistency in successful feedback strengths as a function of problem size for the GD-CV-CIM and momentum-CV-CIM; the step sizes required for converge on instances with small $\kappa$ lead to divergent behavior when $\kappa$ is large.
We expect that, much like classical optimization, problem conditioning is a primary factor dictating stable feedback strengths.
While setting the feedback strength to be too small hampers the convergence of the Adam-CV-CIM, its dynamics tend to be consistent as a function of feedback strength, even across instances of varying $\kappa$.
This is likely because Adam regularizes the momentum updates using the second-moment estimates, and is scale-invariant as a result.
Moreover, the learning rate automatically adapted per variable, allowing Adam to explicitly combat the poor problem scaling when $\kappa$ is large.

Surprisingly, momentum is more unstable than gradient descent, while it is typically known for improving stability for poorly conditioned instances.
This is likely because the CIM relies solely on the double-well potential to enforce the constraints, and they are not otherwise incorporated in the feedback calculation.
As a consequence, the momentum update continues to add up the contribution from the gradients, even if they lead to constraints violation.
This suggest that the feedback term should also take constraints into account, similarly to the approach taken by the AHC-CIM.

The purpose of this experiment is not to suggest that there are no values of $\lambda$ or re-parameterizations of the instances where the performance of the GD-CV-CIM and momentum-CV-CIM are comparable to that of the Adam-CV-CIM.
In fact, a diagonal pre-conditioner that undoes the skew-inducing step in the instance construction would suffice.
However, parameter-tuning can be tedious process that may result in disparate results depending on the expertise and persistence of the user, particularly when differences in problem characteristics prevent extrapolating good parameter settings from one set of instances to another.
While it is unlikely to make these solvers fully parameter-free, a solver that is robust, reliable, and user-friendly will be accommodating of some user-error in the parameter settings.
These qualities are part of the reason for Adam's wide adoption in deep-learning, and we believe that adoption of non-conventional optimization architectures will follow similar desiderata.

\section{Discussion}
In this work, we introduced the CIM and the Gaussian-state model of its dynamics, highlighting clear opportunities to incorporate more sophisticated optimization techniques in the feedback term.
We benchmarked variants with the Adam and momentum updates, replacing the standard feedback term.
We found that both the Adam-CV-CIM and momentum-CV-CIM greatly improved convergence speed over the GD-CV-CIM. 
While they did not consistently improve the probability of reaching specific optimality gap targets within our benchmarks, particularly when all samples from the GD-CV-CIM converged to the global optimum, they generally resulted in a greater diversity of samples.
We also found that the performance of the Adam-CV-CIM was significantly less sensitive to the feedback strength parameter (analogous to step size in classical optimization), resulting in one fewer parameter that needs to be meticulously tuned.

Many of the same benefits that Adam and momentum confer to classical optimization also translate to the CIM.
While this is unsurprising when the exposition of the CIM is couched in the language of optimization, it typically is not, obscuring the connections to classical optimization theory.
Despite the limitations of classical optimization algorithms running on von Neumann computers, they have undoubtedly achieved impressive empirical success, including on difficult, non-convex optimization problems.
Even if non-conventional computing architectures are intended to address the limitations of von Neumann computers, we contend that the lessons learned from classical optimization theory should be re-purposed where appropriate.
For example, some challenges identified by the numerical experiments include: (1) better incorporation of the constraints in the feedback computation, which could potentially be addressed with a Lagrangian-inspired update,  and (2) standardization of problem inputs to ease parameter tuning, which could take inspiration from preconditioning. 

This raises the question of how one should trade off the complexity of $\Phi(\{\tilde{\mu}(t) \}_t)$ and the performance of the solver as a whole.
One can envision increasingly complicated feedback terms, $\Phi(\{\tilde{\mu}(t) \}_t)$, that start to resemble higher-order optimization algorithms at the expense of more computational effort.
We expect a Pareto frontier of feedback terms with different trade-offs between time/energy-per-iteration and the number of iterations required to reach a desired optimization target.
The optimal feedback term for a particular scenario will likely depend on resource constraints, such as wall clock time or energy usage.

In this work, we consider a discrete-time perspective of the dynamical system for ease of simulation; however, the continuous-time dynamical systems perspective has recently had significant impacts on optimization theory \cite{wilson2021lyapunov}.
Exploring the continuous-time models of optimization algorithms will likely result in tighter analogs with the continuous-time models often used to describe non-conventional architecture, which we leave to future work.
\small{
\subsection*{Acknowledgements}
This work was supported by NSF CCF (grant \#1918549), NASA Academic Mission Services (contract NNA16BD14C – funded under SAA2-403506).  
R.B. acknowledges support from the NASA/USRA Feynman Quantum Academy Internship program.}

\bibliographystyle{splncs04}
\bibliography{references}

\end{document}